\pdfoutput=1
\RequirePackage{ifpdf}
\ifpdf 
\documentclass[pdftex]{sigma}
\else
\documentclass{sigma}
\fi

\providecommand{\F}[1]{\mathbb{#1}}
\newcommand{\ZZ}{\mathbb{Z}}
\newcommand{\NN}{\F{N}}
\newcommand{\QQ}{\F{Q}}
\newcommand{\CC}{\F{C}}
\newcommand{\PP}{\F P}
\newcommand{\plus}{\oplus}
\newcommand{\tensor}{\otimes}
\newcommand{\w}{\omega}
\newcommand{\W}{\Omega}
\providecommand{\Cal}[1]{\mathcal{#1}}
\newcommand{\CO}{\Cal O}
\newcommand{\CM}{\Cal M}
\newcommand{\CN}{\Cal N}
\newcommand{\ClC}{\Cal C}
\providecommand{\norm}[1]{\lVert#1\rVert}
\providecommand{\set}[1]{\{#1\}}
\newcommand{\p}{\rho}

\numberwithin{equation}{section}

\newtheorem{Theorem}{Theorem}[section]
 { \theoremstyle{definition}
\newtheorem{Definition}[Theorem]{Definition}
\newtheorem{Exercise}[Theorem]{Exercise}
 }

\begin{document}

\newcommand{\arXivNumber}{1809.05732}

\renewcommand{\thefootnote}{}

\renewcommand{\PaperNumber}{129}

\FirstPageHeading

\ShortArticleName{Aspects of the Topology and Combinatorics of Higgs Bundle Moduli Spaces}

\ArticleName{Aspects of the Topology and Combinatorics\\ of Higgs Bundle Moduli Spaces\footnote{This paper is a~contribution to the Special Issue on Geometry and Physics of Hitchin Systems. The full collection is available at \href{https://www.emis.de/journals/SIGMA/hitchin-systems.html}{https://www.emis.de/journals/SIGMA/hitchin-systems.html}}}

\Author{Steven RAYAN}

\AuthorNameForHeading{S.~Rayan}

\Address{Department of Mathematics \& Statistics, McLean Hall, University of Saskatchewan,\\ Saskatoon, SK, Canada S7N 5E6}
\Email{\href{mailto:rayan@math.usask.ca}{rayan@math.usask.ca}}
\URLaddress{\url{https://www.math.usask.ca/~rayan/}}

\ArticleDates{Received September 23, 2018, in final form December 04, 2018; Published online December 07, 2018}

\Abstract{This survey provides an introduction to basic questions and techniques surrounding the topology of the moduli space of stable Higgs bundles on a Riemann surface. Through examples, we demonstrate how the structure of the cohomology ring of the moduli space leads to interesting questions of a combinatorial nature.}


\Classification{14D20; 46M20; 57N65; 05A19}

\renewcommand{\thefootnote}{\arabic{footnote}}
\setcounter{footnote}{0}

\section{Introduction}

Nonabelian Hodge theory realizes an equivalence between three types of objects in geometry and topology: representations of the fundamental group of a complex projective manifold, flat connections on that manifold, and Higgs bundles on that same manifold. The first type of object is topological, the second records the smooth geometry of the manifold, and the third is holomorphic. The nonabelian Hodge correspondence can be formulated into a diffeomorphism of appropriately-defined moduli spaces of these objects. One of the nice features of working on the ``Higgs'' side is the existence of a Hamiltonian ${\rm U}(1)$-action~-- equivalently, an algebraic $\CC^\star$-action, depending on how exactly one constructs the moduli space. By localization, one can at least in principle compute numerical topological invariants of the Higgs bundle moduli space using this action and then possess, by virtue of nonabelian Hodge theory, these invariants for all three moduli spaces.

While the ${\rm U}(1)$-action provides a place to get started, the localization calculation does not scale easily, with explicit results revealing themselves readily only in low rank, even when we restrict to Riemann surfaces. That being said, the structure of the fixed-point locus hints at interesting combinatorics lurking in the cohomology ring of the moduli space, some of which we see below. The fact that the cohomology ring lies at the centre of a number of conjectures in mirror symmetry~\cite{MR1990670} (some of which have been recently addressed \cite{2017arXiv170708536G, 2017arXiv170706417G}) makes these combinatorial questions even more intriguing.

In this article, we present some basic concepts and examples surrounding the problem of computing topological invariants of Higgs bundle moduli spaces. For simplicity, we restrict to the Betti numbers of the rational cohomology ring. The article is based more or less on a mini-course given by the author at the first ``Workshop on the Geometry and Physics of Higgs Bundles'', held in October 2016 at the University of Illinois at Chicago. The mini-course consisted of three lectures and three problem sessions. The presentation in this article, much as in the mini-course, is somewhat bare bones and involves only traditional Morse--Bott theory. For Higgs bundles this is by now ``old hat'', having been supplanted by a number of refinements or wholly different techniques, including arithmetic harmonic analysis; wall-crossing techniques; and motivic and $p$-adic integration. These techniques have led to explicit results about the cohomology that once seemed quite far away. It is difficult to provide a complete list of references on these developments, although here are some that reflect the evolution of these developments: \cite{MR3293805,2017arXiv170706417G,TH:13,MR3003926,MR2453601, 2017arXiv170704214M,MR2975380,2014arXiv1411.2101M,2017arXiv170504849M,MR3432585}.

The mini-course had been delivered for an audience of mostly beginning graduate students. This survey has been written with similar considerations in mind. We imagine that the reader possessing some basic Riemann surface theory~-- including Jacobians, \v{C}ech cohomology, Serre duality, and the Riemann--Roch theorem for holomorphic vector bundles~-- will get the most from these notes. We have included a few basic exercises to capture some of the spirit of the problem sessions.

\section{Background on Higgs bundles}

Higgs bundles originated within mathematical inquiries into gauge theories in the 1970s and 1980s but can also be understood in a mostly algebraic way. We briefly examine both points of view here, with an aim to understanding roughly the geometric features of the moduli space of Higgs bundles.

\subsection{Gauge theory} From this point forward, $X$ is a smooth compact Riemann surface. For now, the genus $g$ of $X$ is at least $2$. We use the symbols $\mathcal O_X$ and $\omega_X$ for the trivial line bundle and cotangent bundle of $X$, respectively. Higgs bundles originally arose as solutions of the \emph{Hitchin equations} or ``self-duality equations'' on $X$ \cite{MR887284}. These are self-dual, dimensionally-reduced Yang--Mills equations written on a smooth Hermitian bundle of rank $r\geq1$ and degree $0$ on $X$. We will use $E$ for this bundle and $h$ for the metric. The equations take the form
\begin{gather}
F(A)+\phi\wedge\phi^* = 0,\label{Hit1}\\
\overline\partial_A\phi = 0.\label{Hit2}\end{gather} In the equations, $A$ is a connection on the bundle (unitary with regards to $h$), $F$ is its curvature, and $\phi$ is a smooth bundle map from $E$ to $E\otimes\omega_X$, called a \emph{Higgs field}. The equations are trivially satisfied by a flat connection $A$ with $\phi=0$. Equation \eqref{Hit1} says that, whenever $A$ is not flat, its curvature $(1,1)$-form should be expressible in terms of $\phi$ and its Hermitian adjoint. Equation~\eqref{Hit2} says that $\phi$ should be holomorphic with respect to the holomorphic structure on $E$ induced by $A$. The equations can be altered appropriately, involving a constant central curvature term on the right side of~\eqref{Hit1}, in order to accommodate an arbitrary degree $d\in\ZZ$. Throughout, we will assume that $r$ and $d$ are coprime.

Now, assume that $\mathcal E$ is a holomorphic bundle on $X$ together with a holomorphic section $\phi\in H^0(X,\operatorname{End}(\mathcal E)\otimes\omega_X)$. We refer to such a pair as a \emph{Higgs bundle}. One can ask: when does the data $(\mathcal E,\phi)$ arise from a solution to the Hitchin equations? In other words, when does there exist a Hermitian metric $h$ on the underlying smooth bundle and a unitary connection $A$ such that the holomorphic structure on $\mathcal E$ is induced by $(h,A)$ and $(A,\phi)$ is a solution of the Hitchin equations for $h$? The answer is a numerical condition on the pair $(\mathcal E,\phi)$, asking that the following inequality holds: for each subbundle $0\subsetneq\mathcal U\subsetneq\mathcal E$ for which $\phi(\mathcal U)\subseteq\mathcal U\otimes\omega_X$, we must have\begin{gather*}\frac{\deg(\mathcal U)}{\operatorname{rank}(\mathcal U)}<\frac{\deg(\mathcal E)}{\operatorname{rank}(\mathcal E)}.\end{gather*} Such $\mathcal U$ are said to be $\phi$-invariant and the ratio in question is referred to as the \emph{slope} of $\mathcal U$. If the inequality is satisfied for all such $U$, we say that the Higgs bundle $(\mathcal E,\phi)$ is \emph{stable}. (The edge case where equality is permitted, known as semistability, is eliminated by the earlier coprime assumption.)

This correspondence is an example of what are now generally referred to as \emph{Kobayashi--Hitchin correspondences}, relating bundles with special metrics to ones with algebro-geometric restrictions. As an equivalence of moduli spaces, on one side we have the space of solutions $(A,\phi)$ of~\eqref{Hit1} and \eqref{Hit2} for $(E,h)$ taken up to gauge equivalence, which are orbits of the conjugation action of the group of smooth unitary diffeomorphisms of~$E$. This quotient has the structure of a smooth, non-compact manifold. On the Higgs bundle side, we have the space of all stable pairs $(\mathcal E,\phi)$ with underlying smooth bundle $E$ taken up to isomorphism, which is given by the conjugation action of the group of holomorphic automorphisms of $\mathcal E$. This quotient has the structure of a non-singular, quasiprojective variety.

The gauge-theoretic side can be interpreted as an infinite-dimensional hyperk\"ahler quotient, in the sense of \cite{MR877637}. Here, the hyperk\"ahler moment maps are the left side of \eqref{Hit1} and the real and imaginary parts of the left side of \eqref{Hit2}. The quotient inherits a hyperk\"ahler metric, compatible with three quaterionically-commuting complex structures. It is an immediate consequence that the moduli space is Calabi--Yau, although it is not compact. The moduli variety on the other side of the correspondence, which we denote by $\CM_X(r,d)$, can be interpreted as a geometric-invariant theory quotient, with its stability condition given by our notion of ``stable'' above. Indeed, this is exactly the condition required to form a Hausdorff moduli space here.

This correspondence generalizes the earlier one of Narasimhan--Seshadri \cite{MR0184252}, which relates flat bundles to stable holomorphic bundles. At the same time, the Kobayashi--Hitchin correspondence can be viewed as a ``fourth corner'' in nonabelian Hodge theory, extending the equivalence to one between flat connections, representations of $\pi_1(X)$, Higgs bundles, and solutions of Hitchin's equations.

For our purposes (and until we introduce some tools from differential topology in Section \ref{SectU1}), we will lean in an algebro-goemetric direction and concentrate on Higgs bundles and $\CM_X(r,d)$. For a deeper discussion of the gauge theory, including an exploration of recent results concerning the global properties of the hyperk\"ahler metric, we refer the reader to \cite{LF:18} in the same collection of mini-course articles~-- as well as of course Hitchin's original article \cite{MR887284}. Regarding nonabelian Hodge theory in particular, we refer the reader to works of Simpson \cite{MR1159261,MR1179076} and to recent surveys such as \cite{MR3409775, MR3675465}.

One common preference, which is useful for instance when going from Higgs bundles to representations of surface groups, is to fix the determinant of the Higgs bundle, which means taking $\wedge^r\mathcal E$ to be some fixed degree-$d$ line bundle. This takes us from the vector bundle (i.e., ${\rm GL}(r,\CC)$) situation to principal ${\rm SL}(r,\CC)$-Higgs bundles. Accordingly, the Higgs field is taken to be trace-free, which we denote by $\phi\in H^0(X,\operatorname{End}_0(\mathcal E)\otimes\omega_X)$. We will use $\CM^0_X(r,d)$ to denote this moduli space, i.e., that of stable ${\rm SL}(r,\CC)$-Higgs bundles with fixed determinant of degree~$d$.

\subsection{Deformation theory}\label{SectDef} The first piece of topological information to compute about $\CM_X(r,d)$ is its dimension. For this, we can use deformation theory. Let us assume, to begin with, that we are working with ${\rm SL}(r,\CC)$-Higgs bundles. To such a Higgs bundle $(\mathcal E,\phi)$, we can associate a deformation complex determined by the \v{C}ech co-differential $\delta$ on $\mathcal E$ and the Higgs field itself. We can view the Higgs field as a map that acts on Lie-algebra-valued forms by the Lie bracket on the Lie algebra part and by the wedge product on the form part. In our situation, where the Higgs field is a section of $\operatorname{ad}(\mathcal E)\otimes\omega_X\cong\operatorname{End}_0(\mathcal E)\otimes\omega_X$, the fact that $\omega_X\wedge\omega_X=0$ on a curve means that the map $(\wedge\phi)^2$ is always zero and hence is a co-differential for our purposes. (For $X$ of higher dimension, this is one motivation for including an extra condition on Higgs bundles, namely that $\phi$ satisfies $\phi\wedge\phi=0$.)

By analogy with the fact that the tangent space to the moduli space of stable bundles at a~point $E$ is the cohomology $H^1(X,\operatorname{End}_0(\mathcal E))$ of the complex associated to $\delta$, the tangent space to the moduli space at a stable pair $(\mathcal E,\phi)$ is the hypercohomology $\mathbb H^1$ of the double complex associated to the two co-differentials, $\delta$ and $\wedge\phi$~\cite{MR1260109}. By working with the double complex as in~\cite{MR1260109}, we find that $\dim_{\CC}\mathbb H^1$ is a sum of two numbers. The first is the dimension of \begin{gather*}\ker H^1(X,\operatorname{End}_0(\mathcal{E}))\stackrel{\wedge\phi}{\longrightarrow}H^1(X,\operatorname{End}_0(\mathcal{E})\otimes\omega_X),\end{gather*}which is a subspace of the usual tangent space to the moduli space of stable bundles. Here, we only want deformations of the holomorphic structure on the bundle for which $\phi$ is still holomorphic itself. The second number is the dimension of\begin{gather*}\frac{H^0(X,\operatorname{End}_0(\mathcal{E})\otimes\omega_X)}{\operatorname{im}H^0(X,\operatorname{End}_0(\mathcal{E}))\stackrel{\wedge\phi}{\longrightarrow} H^0(X,\operatorname{End}_0(\mathcal{E})\otimes\omega_X)},\end{gather*} which captures deformations of the Higgs field.

It is a consequence of stability that the map
\begin{gather*} \wedge\phi\colon \ H^0(X,\operatorname{End}_0(\mathcal{E})){\longrightarrow}H^0(X,\operatorname{End}_0(\mathcal{E})\otimes\omega_X) \end{gather*} is injective. (See, for instance, \cite[Remark~2.8]{MR3675465}.) It then follows by duality that the map
\begin{gather*} \wedge\phi\colon \ H^1(X,\operatorname{End}_0(\mathcal{E})){\longrightarrow}H^1(X,\operatorname{End}_0(\mathcal{E})\otimes\omega_X)\end{gather*} is surjective.
\begin{Exercise}\label{Exercise1}Show that $\dim_{\CC}\CM^0_X(r,d)=2(r^2-1)(g-1)$.\footnote{\emph{Hint}: Each of the two numbers that must be summed to give $\dim_{\CC}\mathbb H^1$ can be expressed as a difference, owing to the injectivity and surjectivity properties. These differences can be rearranged in such a way that Riemann--Roch can be applied.}
\end{Exercise}

With this in place, it is easy to reason in a number of ways that $\dim_{\CC}\CM_X(r,d)=2r^2(g-1)+2$. The difference between the two dimensions is~$2g$, which is the sum of the dimension of the Jacobian of $X$ and number of linearly independent $1$-forms on $X$~-- the latter accounts for removing the trace from $\phi$.

\subsection{Examples} The Kobayashi--Hitchin correspondence allows us to construct examples of solutions to Hitchin's equations as Higgs bundles, simply by combining a holomorphic bundle with a Higgs field $\phi$ that fails to preserve ``bad'' subbundles with excess slope. One can achieve this by constructing a~Higgs field that leaves \emph{no} proper subbundle invariant whatsoever. In fact, if $\mathcal E=\mathcal L$ is a~holomorphic line bundle on $X$, then any $\phi$ has this property, and so a line bundle with a section $\phi\in H^0(X,\mathcal L\otimes\mathcal L^*\otimes\omega_X)=H^0(X,\w_X)$, which is nothing more than a holomorphic one-form, is an example of a Higgs bundle.

\begin{Exercise}\label{Exercise2} Show that $\CM_X(1,d)$ is homeomorphic to $\mathbb{R}^{2g}\times\big(S^1\big)^{2g}$ and that $\CM^0_X(1,d)$ is just a point.
\end{Exercise}

A more interesting example comes from considering the rank-$2$, degree-$0$ split bundle $\mathcal E\cong\w_X^{1/2}\plus\w_X^{-1/2}$, where $\w_X^{1/2}$ is a choice of holomorphic square root of $\w_X$. (There are $2^{2g}$ such line bundles on $X$.) The anti-diagonal Higgs field\begin{gather*}\phi=\left(\begin{matrix}0 & \alpha\\1 & 0\end{matrix}\right)\end{gather*}preserves neither summand of $\mathcal E$, and so is stable. Here, $1$ is interpreted as the identity endomorphism for $\w_X^{1/2}$. The section $\alpha$ is a quadratic differential on~$X$. Hence, we have injective maps from $H^0\big(X,\w_X^{\otimes2}\big)$ into $\CM^0_X(2,0)$ and $\CM_X(2,0)$. Through the Hitchin equations, the existence of this particular family of Higgs bundles induces a uniformizing metric on~$X$, as in Hitchin's paper~\cite{MR887284}.

\subsection{Hitchin fibration} The principal tool for understanding the structure of $\CM_X(r,d)$ is the \emph{Hitchin map}, which is nothing more than the map that assigns to each Higgs bundle the characteristic polynomial (interpreted correctly) of its Higgs field. We write\begin{gather*}\Theta\colon \ \mathcal M_X(r,d)\longrightarrow\mathcal A_r:=\displaystyle\bigoplus_{i=1}^rH^0\big(X,\omega_X^{\otimes i}\big) \end{gather*} defined by sending the isomorphism class of $(\mathcal E,\phi)$ to the $r$-tuple of coefficients of the characteristic polynomial, each of which is a section of a respective tensor power of $\omega_X$. The codomain~$\mathcal A_r$ is an affine space called the \emph{Hitchin base}. The map $\Theta$ is proper and thus fibres $\CM_X(r,d)$ by compact subvarieties, the \emph{Hitchin fibres}. This properness result was established for the space $\mathcal M_X(2,d)$ by Hitchin \cite{MR887284}. In general, see~\cite{MR1085642}.

This gives us a very coarse idea of how the moduli space ``looks'': it is an affine space populated by compact fibres, the generic ones certainly being smooth. Can we sharpen this? To do so, we take a closer look at the characteristic polynomial of a given $\phi$~-- namely, we want to understand the geometry of its roots. Denote by $|\omega_X|$ the total space of $\omega_X$; by $(x,y(x))$, a~local coordinate on $|\omega_X|$ ($x$ is ``horizontal'' and $y$ is ``vertical''); and by $p$, the bundle projection $\omega_X\to X$. The bundle $\p^*\omega_X$ on $|\omega_X|$ has a natural section $w$ given by $w(x,y(x))=y(x)$, where the output value is seen as living in the copy of the fibre $(\w_X)_x$ attached to itself at $y(x)$ in the pullback bundle. This is the so-called \emph{Seiberg--Witten differential}. These objects allow us to define:

\begin{Definition}\label{Definition1} The \emph{spectral curve} determined by $a=(a_1,\dots,a_r)\in\mathcal A_r$ is the $1$-dimensional subvariety $X_a\subset|\omega_X|$ given by the zero locus of the polynomial\begin{gather*}w^r(y)+a_1(p(y))w^{r-1}(y)+\cdots+a_r(p(y)). \end{gather*}\end{Definition}

For a sufficiently general choice of $a$, $X_a$ is a non-singular curve ramified over $X$ with order $r$. In other words, it is an $r:1$ branched cover and so we have fashioned a new Riemann surface, related to $X$, from data in the Hitchin base $\mathcal A_r$. Now, consider any line bundle $\mathcal L$ on $X_a$. The direct image $p_*\mathcal L$ is a locally-free sheaf of rank $r$ and hence can be identified as the sheaf of sections of a holomorphic bundle $\mathcal E\to X$. The Seiberg--Witten differential, thought of as acting by\begin{align*} w|_{X_a} \colon \ \mathcal L & \longrightarrow \mathcal L\otimes p^*\omega_X, \\ s& \longmapsto s\cdot y \end{align*} on the line bundle, pushes forward to a linear map between the sheaves $\mathcal E$ and $\mathcal E\otimes\omega_X$. In other words, we have constructed a Higgs field $\phi$ for the bundle $\mathcal E$, and so the data of a line bundle on $X_a$ leads to a Higgs bundle on~$X$. In the opposite direction, a Higgs bundle $(\mathcal E,\phi)$ on $X$ determines a tuple $a\in\mathcal A_r$ through the Hitchin map. This tuple generates a spectral curve~$X_a$, which is exactly the spectrum of $\phi$, producing distinct eigenvalues at most points $x\in X$ (corresponding to the $r$ sheets of $X_a$, branching wherever there are repeated eigenvalues). The eigenspaces of~$\phi$, which are generically $1$-dimensional, form a sheaf $\mathcal L$ on~$X_a$, which can be shown to be a line bundle. (See Proposition~4.2(2) in Chapter~2 of~\cite{MR1723384}.)

Essentially, we have that an isomorphism class of holomorphic line bundles $[L]$ on $X_a$ is equivalent to the data of an isomorphism class of Higgs bundles $[(\mathcal E,\phi)]$ on $X$. This is the \emph{spectral correspondence} as developed in \cite{MR998478,MR1397059,MR1397273, MR885778}. It follows from it that the generic fibre~$\Theta^{-1}(a)$ is isomorphic to the Jacobian variety of $X_a$. This Jacobian, however, is not typically the space of degree $0$ line bundles on $X_a$. Rather, their degree is shifted by the ramification. The actual degree $e$ is given by\begin{gather*}e = d-(1-g')+r(1-g), \end{gather*} where $g'$ is the genus of $X_a$. We denote this Jacobian by $\operatorname{Jac}^e(X_a)$~-- it has the same dimension regardless of the value of~$e$.

\begin{Exercise}\label{Exercise3} Derive the above formula for $e$.\footnote{\emph{Hint}: Use the Riemann--Roch theorem in combination with properties of the pushforward operation between two smooth curves, one a branched cover of the other.}
\end{Exercise}

Since the genus $g'$ of $X_a$ is equal to the complex dimension of its Jacobian and since $\Theta^{-1}(a)\cong\operatorname{Jac}^e(X_a)$ for generic $a\in\mathcal A_r$, we can obtain the genus of the generic spectral curve by subtracting the dimension of $\mathcal A_r$ from the dimension of the moduli space. For each power of $\omega_X$, Riemann--Roch reads as
\begin{gather*}h^0\big(X,\omega_X^{\otimes i}\big)-h^0\big(X,\omega_X^{\otimes1-i}\big)=(2i-1)(g-1).\end{gather*}
For each $i>1$, $\omega_X^{\otimes1-i}$ has degree $(1-i)(2g-2)>0$ and so $h^0\big(X,\omega_X^{\otimes1-i}\big)$ vanishes, leaving us with
\begin{gather*}h^0\big(X,\omega_X^{\otimes i}\big)= \begin{cases}g & \text{if }i=1,\\(2i-1)(g-1) & \text{if } i>1.\end{cases}\end{gather*} It follows that\begin{gather*}\dim_\CC\mathcal A_r=r^2(g-1)+1.\end{gather*} We observe that this is exactly half the dimension of $\CM_X(r,d)$, and so $g'$ is also $r^2(g-1)+1$. In the ${\rm SL}(r,\CC)$ case, we subtract $h^0(X,\omega_X)=g$ from the dimension of $\mathcal A_r$ (to remove the trace). We denote this reduced based by~$\mathcal A^0_r$. At the same time, we recall that we subtract~$2g$ from the dimension of $\CM_X(r,d)$ to get that of $\CM^0_X(r,d)$, and so the half-dimensionality of the base persists here. (The spectral curve has the same genus as in the ${\rm GL}(r,\CC)$ case, but the Jacobian is replaced with a smaller-dimensional Prym variety.)

For an example, let us examine the moduli space $\CM^0_X(2,0)$. According to the formulas derived above, it has dimension $6g-6$; the generic spectral curve has genus $g'=4g-3$, which is also the dimension of the base $\mathcal A^0_2$; and the degree of the relevant line bundles on the spectral curve is $e=3g-6$. The Hitchin base is just $H^0\big(X,\omega_X^2\big)$, the space of quadratic differentials, which are the possible determinants of $\phi$. If we take $X$ of genus $g=2$ specifically, then the moduli space is $6$-dimensional, fibering over a $3$-dimensional base, with $X$ covered $2:1$ by a~smooth genus $g'=5$ curve $X_a$ for each generic $a\in H^0\big(X,\omega_X^2\big)$. By the spectral correspondence, line bundles of degree $e=0$ push forward from~$X_a$ to produce stable Higgs bundles on $X$. Recall now the family of Higgs bundles $\mathcal E\cong\w_X^{1/2}\plus\w_X^{-1/2}$ with \begin{gather*}\phi=\left(\begin{matrix}0 & \alpha\\1 & 0\end{matrix}\right)\end{gather*} that live in this moduli space. The map $\Theta$ sends $\phi=\left(\begin{smallmatrix}0 & \alpha\\1 & 0\end{smallmatrix}\right)$ to $-\alpha\in H^0\big(X,\omega_X^2\big)$. These Higgs fields form the \emph{Hitchin section}, intersecting each Hitchin fibre in exactly one point. From the spectral point of view, there is a special line bundle on each $X_a$ that pushes forward to produce an element of this family.

\subsection{Integrable system} The moduli space is a fibration in a different way. If $\mathcal N_X(r,d)$ is the moduli space of stable bundles of rank $r$ and degree $d$ (stable here means that \emph{all} proper subbundles must satisfy the slope condition), then the tangent space $T_{\mathcal E}(\mathcal N_X(r,d))$ at some bundle $\mathcal E$ is\begin{gather*}H^1(X,\operatorname{End}(\mathcal{E}))\stackrel{\text{Serre}}{\cong}H^0(X,\operatorname{End}(\mathcal{E})\otimes\omega_X)^*\end{gather*}and so the cotangent bundle to $\mathcal N_X(r,d)$ is contained inside the moduli space of Higgs bundles. It is important to note there are stable Higgs bundles $(\mathcal E,\phi)$ for which the vector bundle $\mathcal E$ alone is unstable and so the projection $\mathcal M_X(r,d)\longrightarrow\mathcal N_X(r,d)$ is only defined above those Higgs bundles with stable underlying bundle. The symplectic form on $T^*\mathcal N_X(r,d)$ can, however, be canonically extended to one on $\mathcal M_X(r,d)$. (The complex structure on $T^*\mathcal N_X(r,d)$ also extends to $\CM_X(r,d)$ in a compatible way, producing one of the complex structures making up the hyperk\"ahler structure on the moduli space.)

Hitchin proved in \cite{MR885778} that this symplectic structure on $\CM_X(r,d)$ is an algebraically completely integrable Hamiltonian system. In particular, the real and imaginary parts of the components of the Hitchin map $\Theta$ are functionally-independent, Poisson-commuting functions, of which there are sufficiently-many due to the half-dimensionality of~$\mathcal A_r$, thereby providing a complete set of Hamiltonians. The Hitchin fibres are the Liouville tori of the dynamical system. Many known integrable systems can be realized as Hitchin systems, with flows linearizing on the Hitchin fibres. (It is often necessary to allow the genus to be $0$ or $1$ and to puncture $X$ so that $\phi$ develops poles at the punctures. This leads naturally to the parabolic Higgs bundle story, cf.~\cite{MR953833,MR1411302}. See also for~\cite{MR1300764} for Hitchin-type integrable systems in which $\omega_X$ is replaced with other line bundles.)

\section[${\rm U}(1)$-action]{$\boldsymbol{{\rm U}(1)}$-action}\label{SectU1}

The coarse description above is not enough to tell us the global topology of the Hitchin fibration. The fibration is nontrivial, due to the presence of special degenerate fibres, and so the global topology is not simply that of a generic torus fibre (unless $r=1$~-- see Exercise~\ref{Exercise2}). It turns out that only one special fibre really matters: this is the one that we call the ``nilpotent cone'', as we will see below.

To study the topology, we could regard the moduli space as the gauge-theoretic moduli space of solutions to Hitchin's equations, in which case we would employ Morse theory for a~suitable height function. For us, this would be the $L^2$-norm on $\CM_X(r,d)$, which is a multiple of $f(\mathcal E,\phi)=\norm{\phi}^2$ coming from the K\"ahler metric associated to the complex structure extended from $T^*\mathcal N_X(r,d)$ (cf.~\cite{MR2955005, MR1298999,MR887284,MR2882774,MR2425469}). Here, we are concerned with critical points of $f$. If we regard the moduli space as the quasiprojective variety $\CM_X(r,d)$, as we have been doing up until now, then we can employ Bia\l{}ynicki-Birula theory \cite{MR0366940} for an algebraic group action. For us, this is the action\begin{gather*}\lambda\cdot(\mathcal E,\phi)=(\mathcal E,\lambda\cdot\phi)\end{gather*}of $\CC^\star$. Here, we are concerned with fixed points of the action. The two approaches are connected by the following fact: all of the fixed points of the action are fixed points of the compact group ${\rm U}(1)\subset\CC^\star$. Moreover, the height function is a moment map for the ${\rm U}(1)$-action and the fixed points of the ${\rm U}(1)$-action are critical points of $f$ \cite{MR887284}.

We denote by $\CM_X(r,d)^{{\rm U}(1)}$ the fixed points of the ${\rm U}(1)$-action. A stable Higgs bundle $(\mathcal E,\phi)$ belongs to $\CM_X(r,d)^{{\rm U}(1)}$ if and only if there exists a automorphism $A_\lambda$ of $\mathcal E$ so that $A_\lambda\phi A^{-1}_\lambda={\rm e}^{{\rm i}\theta}\phi$ for each $\lambda\in[0,2\pi)$. In other words, a Higgs bundle is fixed if and only there is a change of basis that undoes the action of ${\rm U}(1)$. We would like to have a useful description of these fixed points.

\subsection{Holomorphic chains} Now, suppose that $(\mathcal E,\phi)\in\CM_X(r,d)^{{\rm U}(1)}$. If $A_\lambda$ is the one-parameter family of transformations that corrects for the action, then there is a limiting endomorphism $\Lambda$ that generates this family infinitesimally, i.e., \begin{gather*}\Lambda := {\rm D}_\lambda(A_\lambda) |_{\lambda=0},\end{gather*} where ${\rm D}_\lambda$ is a suitably-defined derivative.

\begin{Exercise}\label{Exercise4}
Show that $ [\Lambda,\phi ] = {\rm i}\phi$.\footnote{\emph{Hint}: Differentiate the fixed-point equation $A_\lambda\phi A^{-1}_\lambda={\rm e}^{{\rm i}\theta}\phi$ using the same derivative.}
\end{Exercise}

It is also possible to argue that, if $\overline\partial_A$ is a $\CC$-linear operator that determines the holomorphic structure on $\mathcal E$, e.g., an operator induced by the unitary connection~$A$ satisfying Hitchin's equations, then $\overline\partial_A$ and $\Lambda$ must be simultaneously diagonalizable. (This comes from the fact that automorphisms $A_\lambda$ act trivially by conjugation on the holomorphic structure, by definition of the ${\rm U}(1)$-action.) It follows that $\mathcal E$ decomposes into eigenspaces of~$\Lambda$.

We will call these eigenspaces $\mathcal B_1,\dots,\mathcal B_n$. Geometrically speaking, these are holomorphic subbundles of $\mathcal E$. Likewise, the eigenvalues of $\Lambda$ are global holomorphic functions on $X$: $s_1,\dots,s_n$, respectively. Now, we take some $\mathcal B_k$ and apply both sides of the identity from Exercise~\ref{Exercise4} to it. We find\begin{gather*}\Lambda(\phi\mathcal B_k) = (s_k+{\rm i})(\phi\mathcal B_k),\end{gather*} where ${\rm i}=\sqrt{-1}$. This indicates that the image of $\mathcal B_k$ under the Higgs field is a subbundle of the eigen-bundle for eigenvalue $s_k+{\rm i}$. In turn, this implies that the eigenspaces are grouped into sequences, with their eigenvalues ordered as $s_k$, $s_k+{\rm i}$, $s_k+2{\rm i}$, and so on. These sequences terminate when the image of an eigen-bundle under $\phi$ is zero (or when we reach the last eigen-bundle). It can be shown that the existence of multiple, disconnected sequences for a fixed point would violate stability, as stable Higgs bundles are irreducible in the sense that they cannot decompose into proper, nonzero Higgs subbundles. Hence, it follows that for a rank-$r$ Higgs bundle $(\mathcal E,\phi)\in\CM_X(r,d)^{{\rm U}(1)}$, there exists a number $n$ such that $\mathcal E=\bigoplus_{k=1}^n\mathcal B_k$ and\begin{gather*}\mathcal B_1\stackrel{\phi_1}{\longrightarrow}\mathcal B_2\tensor\omega_X\stackrel{\phi_2}{\longrightarrow}\cdots\stackrel{\phi_{n-1}}{\longrightarrow}\mathcal B_{n}\otimes(\omega_X)^{\otimes(n-1)}\stackrel{\phi_n}{\longrightarrow}0, \end{gather*} where $\phi_k=\phi|_{\mathcal B_k}$ and $\phi_k$ is not identically zero for $k<n$.

A Higgs bundle admitting a description such as above is referred to as a \emph{holomorphic chain}, cf.~\cite{MR1823573,MR2253535,MR2030379,MR3293805}. Equivalently, such Higgs bundles can be regarded as complex variations of Hodge structure~-- see~\cite{MR1159261}.

This description says that we can write a fixed point in a basis of sections where $\phi$ has the blocks $\phi_i$ arranged sub-diagonally:\begin{gather*}\phi=\left(\begin{matrix}0 & 0 & \cdots & 0 & 0\\\phi_1 & 0 & \cdots & 0 & 0\\0 & \phi_2 & \cdots & 0 & 0\\ & & \ddots & & \\ 0 & 0 & \cdots & \phi_{n-1} & 0\end{matrix}\right).\end{gather*} Such a matrix is nilpotent and so every fixed point belongs to the Hitchin fibre $\Theta^{-1}(0)$, which is what we refer to as the \emph{nilpotent cone}. In general, not every point in the nilpotent cone is fixed: only those admitting a strict block sub-diagonal (or super-diagonal) description are fixed.

\begin{Exercise}\label{Exercise5} Show that a Higgs bundle $(\mathcal E,\phi)$ with strict block sub-diagonal Higgs field is necessarily fixed under the ${\rm U}(1)$-action.
\end{Exercise}

If $(\mathcal E,\phi)\in\CM_X(r,d)^{{\rm U}(1)}$, then there is a well-defined $n$-tuple $(r_1,\dots,r_n)$ that encodes the ranks of the $\mathcal B_k$ subbundles~-- this is the \emph{rank vector} of the fixed point.

\subsection{Localization} The key result for us is that the total space of the Hitchin fibration $\CM_X(r,d)$ deformation retracts, via the gradient flow of the moment map of the ${\rm U}(1)$-action, onto $\Theta^{-1}(0)$ \cite{TH:98}. In terms of invariants, the cohomology ring localizes to the fixed-point locus inside $\Theta^{-1}(0)$. The Poincar\'e series ${\rm P}[\CM_X(r,d)]$ that generates the Betti numbers of the rational cohomology $H^\bullet(\CM_X(r,d),\QQ)$ will be a weighted sum of the Poincar\'e series ${\rm P}[\ClC_i]$ of the connected components $\ClC_i$, $i\in I$, of the fixed-point locus. Also, let\begin{gather*}\iota\colon \ \CM_X(r,d)^{{\rm U}(1)}\to\NN\end{gather*}be the function that assigns to each fixed point the number of negative eigenvalues of the Hessian of $f$ at that point, where $f$ is again the moment map. This function $\iota$ is constant on each $\ClC_i$ as per Lemma~9.2 in~\cite{MR1990670} and so the natural number $\iota(\ClC_i)$ is well-defined. It is also worth noting that the rank vector $(r_1,\dots,r_n)$ is constant on connected components of the fixed-point locus, as are the degrees of the~$\mathcal B_k$'s.

Computing $\iota$ will be an important ingredient in the weighted sum that yields ${\rm P}[\CM_X(r,d)]$. Thinking of $\iota$ as the dimension of the ``downward'' subbundle of the normal bundle to $\CM_X\!(r{,}d)^{{\rm U}(1)}\!$ at a fixed point, we can obtain the value of $\iota$ by taking a deeper look at the deformation theory from Section~\ref{SectDef} in the case of a fixed point (cf.\ Section~2.1 of~\cite{MR3610261}). When $(\mathcal E,\phi)$ is fixed, so that a decomposition into an ordered sequence of subbundles $\mathcal B_k$ exists, the action of $\phi$ is with weight $1$ with respect to this sequence, i.e.,\begin{gather*}\phi_k\colon \ \mathcal B_k\longrightarrow\mathcal B_{k+1}\otimes\omega_X.\end{gather*} In other words, elements\begin{gather*}\theta\in\frac{H^0(X,\operatorname{End}_0(\mathcal{E})\otimes\omega_X)}{\operatorname{im}H^0(X,\operatorname{End}_0(\mathcal{E})) \stackrel{\wedge\phi}{\longrightarrow}H^0(X,\operatorname{End}_0(\mathcal{E})\otimes\omega_X)}\end{gather*}
that act with weight $\ell=1$ with respect to the sequence form part of the tangent space at $(\mathcal E,\phi)$ to $\CM_X(r,d)^{{\rm U}(1)}$. The other part comes from the elements\begin{gather*}\beta\in\ker H^1(X,\operatorname{End}_0(\mathcal{E}))\stackrel{\wedge\phi}{\longrightarrow}H^1(X,\operatorname{End}_0(\mathcal{E})\otimes\omega_X)\end{gather*} that act with weight $m=0$ on the sequence, preserving the holomorphic structure of each $\mathcal B_k$. (Since the Higgs field is nilpotent, we can use $\operatorname{End}_0$ here regardless of whether the group is ${\rm GL}(r,\CC)$ or ${\rm SL}(r,\CC)$.) The downward flow comes from weights $(\ell,m)$ with $\ell\geq2$ and $m\geq1$. These weights shorten the holomorphic chain until its length is $n=1$ and the Higgs field is zero, taking us to the ``bottom'' of the nilpotent cone. Out of this comes something computational: $\iota(\ClC_i)$ is the sum of the (real) dimensions of the respective $\ell\geq2$ and $m\geq1$ subspaces of the tangent space.

With all of this in place, the localization identity takes the precise form:

\begin{Theorem}[Hitchin \cite{MR887284}] ${\rm P}[\CM_X(r,d)](t) = \sum\limits_{i\in I}t^{\iota(\ClC_i)}{\rm P}[\ClC_i](t)$.\end{Theorem}

Were the moduli space compact, we would have ${\rm P}[\ClC_i](t)=1$ for each $i\in I$, as in standard Morse theory, and so the Poincar\'e series would reduce to $\sum\limits_{i\in I}t^{\iota(\ClC_i)}$. However, in our case the~$\ClC_i$ are generally positive-dimensional with nontrivial contributions to the cohomology ring. For example, the downward flow of $f$ terminates at the points with $\iota=0$, which is also where $\norm{\phi}^2=0$. These global minimizers are precisely the stable Higgs bundles of the form $(\mathcal E,0)$, which is the set of fixed points with rank vector~$(r)$. This component is in fact the moduli space of stable bundles, $\CN_X(r,d)$, which is positive-dimensional for $g\geq1$. For example, if we consider the ${\rm SL}(2,\CC)$ case with fixed determinant of odd degree~$d$, then the Poincar\'e polynomial of this component is known by \cite{MR702806, MR0364254} to be\begin{gather*}{\rm P}\big[\CN^0_X(2,d)\big](t)=\frac{\big(1+t^3\big)^{2g}-t^{2g}(1+t)^{2g}}{\big(1-t^2\big)\big(1-t^4\big)}.\end{gather*} Like the presentation here, \cite{MR702806} also takes a Morse-theoretic approach. The Poincar\'e series of $\CN_X(r,d)$ factors as the product of ${\rm P}\big[\CN^0_X(2,d)\big](t)$ and that of the Jacobian of $X$ (cf.~\cite{MR702806}), and so we have\begin{gather*}{\rm P}[\CN_X(2,d)](t)=(1+t)^{2g}\frac{\big(1+t^3\big)^{2g}-t^{2g}(1+t)^{2g}}{\big(1-t^2\big)\big(1-t^4\big)}.\end{gather*}

The connected components with higher values of~$\iota$, for which less is immediately known, are an obstruction to determining ${\rm P}[\CM_X(r,d)]$ in high rank, although much recent progress has been achieved via other means as highlighted in the introduction. To shed some light on the difficulty, we recognize that the fixed points can be thought of as representations of $A$-type quivers, with lengths and labels determined by partitions of $r$ and~$d$: \begin{gather*}\bullet_{r_1,d_1}\longrightarrow\bullet_{r_2,d_2}\longrightarrow\cdots\longrightarrow\bullet_{r_n,d_n}.\end{gather*}However, we are not looking at representations in the usual category of vector spaces; rather, we are in the category of bundles on a fixed curve $X$ with $\omega_X$-twisted morphisms. These representations are also known as \emph{quiver bundles}, cf.~\cite{PBG:95,MR2108248,MR3610261,MR3812708,MR2120597}. The moduli space of stable bundles is the solution to the simplest version of this problem, where the quiver has a~single node:\begin{gather*}\bullet_{r,d}.\end{gather*}

Nevertheless, we wish to exhibit a couple of sample calculations in low rank where we can determine this polynomial completely.

\section{Calculations}

\subsection[Rank $r=1$]{Rank $\boldsymbol{r=1}$}

We start off with the simplest possible example, just to have an instance where the answer is readily seen to be correct. The only partition of $r=1$ is the rank vector $(1)$. The entire fibre $\Theta^{-1}(0)$ of $\mathcal M_X(1,d)$, which is the submanifold $\{(\mathcal L, 0)\colon \mathcal{L}\in\operatorname{Jac}^d(X)\}$, is fixed by the ${\rm U}(1)$-action. Hence, there is a single connected component of the fixed-point locus and the number $\iota$ is $0$~-- there are no further components to which to flow down. It follows that\begin{gather*}{\rm P}[\mathcal M_X(1,d)](t) = {\rm P}\big[\operatorname{Jac}^d(X)\big](t) = (1+t)^{2g}, \end{gather*} agreeing exactly with Exercise~\ref{Exercise2}. (Of course, for $\mathcal M^0_X(1,d)$ the moduli space is just a point and the result is even more trivial.)

\subsection[Rank $r=2$]{Rank $\boldsymbol{r=2}$}

\looseness=-1 Now, we look at $\mathcal M_X(2,d)$ for some odd $d$. For convenience, we take $d=1$. Here, we mostly follow Hitchin in~\cite{MR887284}, although there are a few notable differences: we do the ${\rm GL}(2,\CC)$ case rather than ${\rm SL}(2,\CC)$ and our calculation of $\iota$ will use the approach outlined in the preceding section.

The elements of the fixed point set are of two types, $(2)$ and $(1,1)$. Those with rank vector~$(2)$ correspond to the moduli space of stable bundles on $X$, as mentioned earlier. These are the fixed points with $\iota=0$, as per the previous section. Therefore, the contribution to the Poincar\'e series is\begin{gather*}t^0(1+t)^{2g}\frac{\big(1+t^3\big)^{2g}-t^{2g}(1+t)^{2g}}{\big(1-t^2\big)\big(1-t^4\big)}.\end{gather*}

Now, each holomorphic chain of type $(1,1)$ consists of two line bundles $\mathcal B_1$ and $\mathcal B_2$ together with a map $\phi_1\colon \mathcal B_1\to\mathcal B_2\otimes\omega_X$. Let $b=\deg\mathcal B_1$, in which case $\deg\mathcal B_2=1-b$. Note that $\mathcal B_2$ is annihilated by the overall Higgs field, and so we must have $1-b$ strictly less than the slope of $\mathcal E=\mathcal B_1\oplus\mathcal B_2$. Hence, $b\geq1$. On the other hand, if $\phi_1=0$, then $\mathcal B_1$ would be invariant, which violates stability as $b$ would exceed the slope of $\mathcal E$. Having $\phi\neq0$ requires that\begin{gather*}\deg(\mathcal B_1^*\otimes\mathcal B_2\otimes\omega_X)=2g-2b-1\end{gather*}is nonnegative. Taking these together, we have $1\leq b\leq g-1$.

Certainly, two choices of $\mathcal B_1$ with different degrees cannot lie in the same connected component of $\CM_X(2,1)^{{\rm U}(1)}$. Therefore, let us fix a value of $b$ in the range above. The data is thus a triple of a line bundle in $\operatorname{Jac}^{b}(X)$, another in $\operatorname{Jac}^{1-b}(X)$, and a map in $H^0(X,\mathcal B_1^*\otimes\mathcal B_2\otimes\omega_X)$. The dimension of the third space depends on $\mathcal B_1$ and $\mathcal B_2$. To clarify this, suppose~$\mathcal B_2$ is fixed. Instead of keeping track of $\mathcal B_1$, we can instead deal with $\mathcal D=\mathcal B_1^*\otimes\mathcal B_2\otimes\omega_X$. The choice of $\mathcal B_1$ determines~$\mathcal D$ and vice-versa. The relevant data is now the pair $(\mathcal D,\phi_1)$ in which~$\mathcal D$ is a line bundle of degree $-2b+2g-1$ and $\phi_1$ is a holomorphic section of this line bundle. Since $\phi_1$ is not identically zero, this data determines an effective divisor of degree $-2b+2g-1$ on~$X$, which is an element of the $(-2b+2g-1)$-fold symmetric product of~$X$ with itself: ${\rm S}^{-2b+2g-1}(X)$. Notice that for $g\geq2$ and $1\leq b\leq g-1$, the order of this product is always positive~-- in other words, we are considering divisors of at least $1$ point. An element of this symmetric product determines a line bundle~$\mathcal D$ together with a nonzero section~$\phi_1$ vanishing on the divisor. This section is determined only up to scale, i.e., $\phi_1\in\mathbb P H^0(X,\mathcal D)$. However, since we are working inside the moduli space $\CM_X(2,1)$, we are only considering holomorphic chains up to equivalence by automorphisms of $\mathcal E=\mathcal B_1\oplus\mathcal B_2$ that preserve the structure of a $(1,1)$ chain. In other words, we are free to use the action of $\mathbb C^*\times\mathbb C^*\subset\operatorname{Aut}(\mathcal E)$ to put a given chain into a representative form. We can use either~$\mathbb C^*$ to identify any two $\phi_1$'s that differ only by scale, and so the projective representatives given by the divisor coincide exactly with the equivalence classes of pairs $(\mathcal D,\phi_1)$ in the moduli space.

Hence, $\CM_X(r,d)^{{\rm U}(1)}$ has $g$ connected components: the moduli space of stable bundles together with $g-1$ components coming from fixed points with rank vector $(1,1)$. By the argument above, components of the latter type are indexed by $b$ in $1\leq b\leq g-1$ and each component is a bundle over ${\rm S}^{-2b+2g-1}(X)$ with fibre $\operatorname{Jac}^{1-b}(X)$, where the Jacobian accounts for the choice of $\mathcal B_2$. For each $b$ we need the Poincar\'e series of the respective $(-2b+2g-1)$-fold symmetric product of $X$. These generating functions are due to Macdonald~\cite{MR0151460}. Specifically, the Poincar\'e polynomial, in~$t$, of ${\rm S}^nX$ is the coefficient of $s^n$ in the Taylor--Maclaurin series expansion of\begin{gather*}\frac{(1+st)^{2g}}{(1-s)\big(1-st^2\big)}.\end{gather*}

\looseness=-1 Now, regarding the indices $\iota$ for the type $(1,1)$ components, we note that the only element $\theta$ in\begin{gather*}\frac{H^0(X,\operatorname{End}_0(\mathcal{E})\otimes\omega_X)}{\operatorname{im}H^0(X,\operatorname{End}_0(\mathcal{E})) \stackrel{\wedge\phi}{\longrightarrow}H^0(X,\operatorname{End}_0(\mathcal{E})\otimes\omega_X)}\end{gather*} acting with weight $2$ or higher on the sequence $(\mathcal B_1,\mathcal B_2)$ is $\theta=0$, as there are only two bundles in the sequence. Hence, we need only account for elements $\beta$ of weight at least $1$ in\begin{gather*}\ker H^1(X,\operatorname{End}_0(\mathcal{E}))\stackrel{\wedge\phi}{\longrightarrow}H^1(X,\operatorname{End}_0(\mathcal{E})\otimes\omega_X).\end{gather*}For the same reasons, there are no elements of weight~$2$ or higher, and so we seek the elements of weight exactly $1$. Before the action of $\wedge\phi$, the weight $1$ elements form $H^1(X,\mathcal B_1^*\otimes\mathcal B_2)$. The map~$\wedge\phi$ sends these to weight $2$ elements in $H^1(X,\mathcal B_1^*\otimes\mathcal B_2\otimes\omega_X)$. Since the only weight $2$ element is the zero element, we have that all weight $1$ elements are in the kernel of $\wedge\phi$. Our calculation of $\iota$ thereby reduces to the real dimension of $H^1(X,\mathcal B_1^*\otimes\mathcal B_2)$. Since $\deg(\mathcal B_1^*\otimes\mathcal B_2)=1-2b<0$, we have that $H^0(X,\mathcal B_1^*\otimes\mathcal B_2)$ vanishes. Then, by Riemann--Roch we have\begin{gather*}\iota(\mathcal E,\phi) = 4b-4+2g.\end{gather*}

Taking all of this together, we get that the Poincar\'e series of $\CM_X(2,1)$ is\begin{gather*}{\rm P}[\CM_X(2,1)][t]\\
\qquad{} = (1+t)^{2g}\left(\frac{\big(1+t^3\big)^{2g}-t^{2g}(1+t)^{2g}}{\big(1-t^2\big)\big(1-t^4\big)}+\sum_{b=1}^{g-1}t^{4b-4+2g}{\rm P}\big[{\rm S}^{-2b+2g-1}(X)\big](t)\right),\end{gather*}where the Poincar\'e polnyomials for the symmetric products come from Macdonald's function.

\begin{Exercise}\label{Exercise6}
Using the results above, check that when $g=2$, we have that\begin{gather*}{\rm P}[\CM_X(2,1)][t]=(1+t)^4\big(1+t^2+4t^3+2t^4+4t^5+2t^6\big).\end{gather*}
\end{Exercise}

\begin{Exercise}\label{Exercise7}Using the results above, check that when $g=3$, we have that
\begin{gather*}{\rm P}[\CM_X(2,1)][t]=(1+t)^6\big(1 + t^2 + 6 t^3 + 2 t^4 + 6 t^5 + 17 t^6 + 12 t^7 + 18 t^8 + 32 t^9\\
\hphantom{{\rm P}[\CM_X(2,1)][t]=}{} + 18 t^{10} + 12 t^{11} + 3 t^{12}\big).\end{gather*}
\end{Exercise}

Notice that the Poincar\'e polynomials above are not palindromes, even though the moduli spaces are smooth. This is of no concern, given that the moduli spaces are non-compact. For example, in $g=3$ the unequal Betti numbers in degrees $0$ and $18$ tell us that, while $\CM_X(2,1)$ is topologically connected ($b_0=1$), the space has a number of irreducible or ``algebraic'' components ($b_{18}=3$ of them). It is also worth noting that the highest power of $t$ in each case is equal to $2r^2(g-1)+2$, which is the \emph{real} dimension of the fibre of the Hitchin map. This is consistent with the fact that the Hitchin base is contractible and the nontrivial topology lies in $\Theta^{-1}(0)$.

A reasonable question is whether ${\rm P}[\CM_X(2,1)](t)/(1+t)^{2g}$ is the Poincar\'e series of $\CM^0_X(2,1)$, the ${\rm SL}(2,\CC)$ moduli space. In general, this is not the case. Rather, the quotient is the generating function for the Betti numbers of the Langlands dual moduli space; that is, the ${\rm PGL}(2,\CC)$ moduli space. The issue is that there is a nontrivial action of the finite group $\Gamma$ of $2$-torsion line bundles~-- the line bundles $\mathcal P$ with $\mathcal P^{\otimes2}=\mathcal O_X$~-- on $\CM^0_X(2,1)$. As a result, there is a variant cohomology and an invariant cohomology with regards to this action. The quotient of $\CM^0_X(2,1)$ by $\Gamma$, which has order $2^{2g}$, is the ${\rm PGL}(2,\CC)$ moduli space. It possesses only the invariant cohomology, whose ranks are given by the coefficients of ${\rm P}[\CM_X(2,1)](t)/(1+t)^{2g}$. For genus $g=2$, this invariant part is\begin{gather*}1+t^2+4t^3+2t^4+4t^5+2t^6,\end{gather*} as in the exercise above. In contrast, the Poincar\'e series of $\CM^0_X(2,1)$ for $g=2$ is\begin{gather*}1+t^2+4t^3+2t^4+34t^5+2t^6\end{gather*}as computed by Hitchin in \cite{MR887284}. Here, we can see the $\Gamma$-variant cohomology concentrating in the degree $5$ part of the cohomology ring. In terms of the calculations, the main difference relative to above is that we are fixing the determinant of $\mathcal E$ to be some fixed line bundle~$\mathcal V$, from which~$\mathcal B_1$ and $\mathcal B_2$ are related by $\mathcal B_2=\mathcal B_1^*\otimes\mathcal V$. Then, to bring in divisors, we need to define a line bundle $\mathcal D=(\mathcal B_1^*)^2\otimes\mathcal V\otimes\mathcal\omega_X$. It follows that instead of symmetric products of $X$, we get $2^{2g}$-fold covers of symmetric products, with fibres consisting of the line bundles~$\mathcal B_1$ whose squares are isomorphic to one another. Here, we see the action of $\Gamma$ working itself into the cohomology.

For further information on the variant versus invariant cohomology, we refer the reader to~\smash{\cite{MR2453601, MR1990670}}. It is also perhaps crucial to point out that the appearance of Langlands duality here is neither superficial nor a red herring. For how Langlands duality manifests in Higgs bundle moduli spaces~-- and how it relates to mirror symmetry~-- we refer the reader to the same reference in addition to \cite{MR2537083,MR2957305,MR2306566}.

The next logical step would be to try our hand at rank $3$. The calculation using Morse theory is noticeably more difficult, because of fixed points with rank vectors $(1,2)$ and $(2,1)$. The type~$(3)$ case remains the moduli space of bundles, whose topological contribution we already know as per above, while the type $(1,1,1)$ fixed points involve symmetric products of~$X$ in an analogous way to the preceding calculations. For $(1,2)$ and $(2,1)$, the data of the fixed point can be converted into a pair $(\mathcal D,\theta)$ in which $\mathcal D$ is a rank $2$ bundle related to the bundles in the chain and $\theta$ is a section of $\mathcal D$. The issue now is to understand the moduli space of such pairs on~$X$. Gothen's approach~\cite{MR1298999} uses Thaddeus' strategy of varying a stability parameter and then constructing the moduli space in steps by keeping track of birational transformations as the parameter is deformed~\cite{MR1273268}. This stability parameter, which is natural in quiver bundle moduli problems, originates in~\cite{MR1124279}. The rank $4$ Poincar\'e series was computed in~\cite{MR3293805} using a~method that is formally similar to the Morse localization above, but which is rooted in motivic considerations. Notably, the $(2,2)$ case had not submitted readily to the variation-of-stability approach, but was resolved via the motivic approach.

We can also ask about the exact structure of the ring $H^\bullet(\mathcal M_X(r,d),\QQ)$ itself. For $r=2$, the generators and relations are worked out in \cite{MR1949162,MR2044052,MR1887889}. For the status of this in higher rank, we refer the reader to~\cite{MR2786591,2018arXiv180810311C}. For examples of Betti numbers over other fields, we refer the reader to \cite{MR3766464} where the $\ZZ_2$ Betti numbers are calculated for rank $2$ Higgs bundles

\section{Combinatorial questions}

In the Morse-theoretic calculations of the preceding section, the degree $d$ of the Higgs bundles enters the calculations explicitly when we work with stable holomorphic chains. However, nonabelian Hodge theory forces the Betti numbers of $\CM_X(r,d)$ to be independent of $d\in\ZZ$, at least when $d$ is coprime to $r$ as we have been assuming all along. This is due to the fact that the Poincar\'e series of the ${\rm GL}(r,\CC)$ character variety of $X$ is insensitive to $d$, where $d$ is used to define twisted representations of $\pi_1(X)$~\cite{MR2453601}. This is combinatorially interesting because there is nothing at first glance to say that corresponding connected components of $\CM_X(r,d)^{{\rm U}(1)}$ have identical Poincar\'e polynomials~-- or even that there are the same number of components.

The $d$-independence of Betti numbers leads to a number of combinatorial observations. We offer a small sample. For our purposes, these are easier to see if we permit $X$ to have genus $g=0$ and if permit Higgs fields twisted by a line bundle other than $\omega_X$. Namely, we wish to consider ``twisted'' Higgs bundles of the form $(\mathcal E,\phi)$ with $\mathcal E$ a vector bundle on the projective line $\mathbb P^1$ and\begin{gather*}\phi\colon \ \mathcal E\longrightarrow\mathcal E\otimes\mathcal O(q),\end{gather*}where $\mathcal O(q)$ is the unique (up to isomorphism) line bundle on $\mathbb P^1$ of degree $q>0$. (The cotangent bundle $\omega_{\PP^1}$ is unsuitable here, as we will then have $q=-2$ and all Higgs bundles of rank $r>1$ and coprime degree $d$ will be unstable.) These Higgs bundles do not rise in the same natural way in gauge theory, but they are nonetheless useful as a~test case here. In particular, these moduli spaces, which are constructed using slope stability in exactly the same way as $\CM_X(r,d)$, have the same natural ${\rm U}(1)$-action~\cite{MR3610261,MR3812708}.

Interestingly, this moduli space does not fit in a natural way into nonabelian Hodge theory~-- one would have to puncture $\PP^1$ along a divisor $D$ and then regard $\phi$ as being valued in $\CO(q)=\omega_X\otimes\CO(D)$ with poles along $D$, with certain conditions on the residues of $\phi$ at the poles \cite{MR2004129, MR1040197}. However, this changes the topology of the moduli space in a significant way and reintroduces the bundle moduli (as we are now keeping track of data in the fibres of $\mathcal E$ at the poles). Keeping our definition the way it is, i.e., holomorphic bundles with holomorphic $\CO(q)$-valued Higgs fields, there is no immediate relationship to a character variety and, as such, no obvious reason for degree independence of the Betti numbers. Yet, it seems to hold in direct calculations of the Betti numbers in low rank, as in~\cite{MR2004129, MR2975380,MR3158239,MR1040197}.

In this setting, because of the relative lack of vector bundle moduli, we attain fairly clear combinatorial descriptions for certain Betti numbers. It is possible for this moduli space to establish via Morse theory that the top Betti number~-- that is, the coefficient of the highest power of $t$ appearing in the Poincar\'e series~-- is precisely the number of connected components of the fixed-point locus coming from fixed points of type $(1,\dots,1)$. This can be shown in turn to be the number of solutions $(d_1,\dots,d_r)\in\ZZ^r$ to the equation\begin{gather*}d_1+\cdots+d_r=d\end{gather*}subject to $d_i-d_{i-1}\leq q$ and, if $r>1$, $(d_j+\cdots+d_r)/(r-j+1)<d/r$ for all $2\leq j\leq r$. Because the $d_j$'s are degrees of line bundles, they are permitted to be negative, and so the equation $d_1+\cdots+d_r=d$ alone is an \emph{unbounded} integer partition problem. The problem becomes well-posed precisely because of stability.

The degree independence of the Betti numbers would, as a corollary, make the solution of this partition problem independent of $d$, again assuming coprimality with regards to $r$. If we fix, say, $q=1$ and then compute the solutions of the above partition problem for increasing $r$, we find the following sequence regardless of which (coprime) $d$ we choose:\begin{gather*}1, 1, 1, 2, 5, 13, 35, 100, 300, 925, 2915, 9386, 30771, 102347, 344705,\dots.\end{gather*} Interestingly, this sequence appears in the OEIS database as \href{https://oeis.org/A131868}{A131868} \cite{OEISA131868}. The entry gives the following function that yields these numbers for each $r$:\begin{gather*} \W(r)=\frac{1}{2r^2}\sum_{e|r}\mu(r/e)\left(\begin{matrix}2e\\ e\end{matrix}\right)(-1)^{e+1},\end{gather*} where $\mu$ is the M\"obius function. By examining type $(1,\dots,1)$ fixed points for other values of $q$ and experimenting with the function $\Omega$, it is not hard to make an educated guess as to a more general version of this function for any~$q$:\begin{gather*} \W(r,q)=\frac{1}{(q+1)r^2}\sum_{e|r}\mu(r/e)\left(\begin{matrix}(q+1)e\\e\end{matrix}\right)(-1)^{qe+1}.\end{gather*} That this is the correct function for all $r>0$, $q>0$ for our counting problem is actually established by Reineke in~\cite{MR2801406}. This also establishes the $d$ independence.

The OEIS entry provides a combinatorial interpretation for the top Betti numbers of the $q=1$ moduli spaces that, while similar in spirit, is not exactly the same as the ours: $r\cdot\Omega(r,1)$ is the number of size $r$ subsets of $\set{1,\dots,2r-1}$ that sum to $1$ modulo $r$. Right away, the degree independence means that we can replace $1\,{\rm mod}\,r$ in this problem with $d\,{\rm mod}\,r$ without changing the solutions. This problem falls into a set of related combinatorial problems studied by Erd\"os--Ginzburg--Ziv~\cite{MR3618568}; in some of these, it is known that one can shift the interval $\set{1,\dots,2r-1}$ freely to any consecutive $2r-1$ numbers (cf.\ the related entry, \href{https://oeis.org/A145855}{A145855} \cite{OEISA145855}). That being said, the partition problem of type $(1,\dots,1)$ fixed points is one in which the \emph{differences} between consecutive parts of the partition are bounded, rather than overall interval in which the parts are allowed to lie.

We can also examine the Poincar\'e series itself as $r$ and $q$ grow. With $r$ fixed and $q$ allowed to grow indefinitely, the Poincar\'e series can be seen to tend to that of the classifying space of the gauge group of the underlying smooth bundle. If we fix $q$ and drive $r$ to larger values~-- or drive both to infinity~-- the series tends to\begin{gather*} 1+t^2+3t^4+5t^6+10t^8+16t^{10}+29t^{12}+45t^{14}+75t^{16}+115t^{18}+\cdots,\end{gather*} whose coefficients are captured in \href{https://oeis.org/A000990}{A000990} \cite{OEISA000990}. If the equivalence of counting problems is correct, this would say that the coefficient of $t^{2n}$ is the number of plane partitions of $n$ with at most $2$ rows. This is especially interesting because it provides a combinatorial interpretation for each Betti number individually, while Morse theory builds each coefficient from potentially many separate combinatorial problems as data from different components of the fixed-point locus contribute to the same coefficient.

Finally, it is worth commenting that in all of these cases~-- the ordinary Higgs bundles of the preceding sections and the twisted ones on $\PP^1$ here~-- that the lack of palindromy in the Poincar\'e series is skewed in such a way that the largest Betti number lies to the ``right'' of the middle coefficients, i.e., between the middle and the top Betti number. This phenomenon is studied in \cite{MR3364745} in the context of non-compact, hyperk\"ahler semiprojective moduli spaces $\mathcal X$. Here, ``semiprojective'' refers to the property of the having an algebraic $\CC^\star$-action with projective fixed-point set with the limit $\lim\limits_{\lambda\to0}\lambda x$ existing for all $x\in\mathcal X$. The fact that this persists for the twisted Higgs bundle moduli spaces on $\PP^1$, which are semiprojective but have no hyperk\"ahler structure, suggests there could be a combinatorial explanation for the phenomenon, independent of the geometry.

In general, we see that for Higgs-bundle-type moduli spaces there is a complicated dance between geometry and combinatorics playing out within the cohomology ring, with geometric phenomena forcing combinatorial identities to emerge and with combinatorial identities expressing themselves geometrically in surprising ways. Throughout, topology is the conduit.

\subsection*{Acknowledgements} I thank Laura Schaposnik for organizing the series of workshops in which the mini-course took place, and both her and Lara Anderson for encouraging the preparation of this survey. With regards to the workshops, I acknowledge support from UIC NSF RTG Grant DMS-1246844, the UIC Start-Up Fund of L.~Schaposnik, and the grants NSF DMS 1107452, 1107263, 1107367 RNMS: GEometric structures And Representation varieties (the GEAR Network). I am grateful to Marina Logares, who gave a mini-course in parallel to mine, for insightful discussions as well as to Laura Fredrickson for useful comments on the manuscript during its preparation. I~thank the referees for helpful remarks and corrections that led to the final version of this article.

\pdfbookmark[1]{References}{ref}
\LastPageEnding

\end{document}